\documentclass[11pt]{amsart}
\usepackage{amssymb,latexsym,amsmath}

\newtheorem{definition}{Definition}[section]
\newtheorem{proposition}{Proposition}[section]
\newtheorem{lemma}{Lemma}[section]
\newtheorem{theorem}{Theorem}[section]
\newtheorem{corollary}{Corollary}[section]

\begin{document}

\title{Integration over the Pauli quantum group}
\author{Teodor Banica}
\address{T.B.: Department of Mathematics, Paul Sabatier University, 118 route de Narbonne, 31062 Toulouse, France. {\tt banica@picard.ups-tlse.fr}}
\author{Beno\^\i{}t Collins}
\address{B.C.: Department of Mathematics, Claude Bernard University, 43 bd du 11 novembre 1918, 69622 Villeurbanne, France and University of Ottawa, 585 King Edward, Ottawa, ON K1N 6N5, Canada. {\tt collins@math.univ-lyon1.fr}}

\begin{abstract}
We prove that the Pauli representation of the quantum permutation algebra $A_s(4)$ is faithful. This provides the second known model for a free quantum algebra. We use this model for performing some computations, with the main result that at the level of laws of diagonal coordinates, the Lebesgue measure appears between the Dirac mass and the free Poisson law.
\end{abstract}

\maketitle

\section*{Introduction}

The notion of free quantum group appeared in Wang's papers
\cite{wa1}, \cite{wa}. The idea is that given a compact group $G\subset
U_n$, the matrix coordinates $u_{ij}\in C(G)$ commute with
each other, and satisfy certain relations $R$. One can define then the universal algebra $A$ generated by abstract variables $u_{ij}$, subject to the relations $R$. In
certain situations we get a Hopf algebra in the sense of Woronowicz \cite{wo}, and we have the heuristic formula $A=C(G^+)$. Here $G^+$
is a compact quantum group, called free version of $G$. The free version is not unique, because it depends on $R$.

This construction is not axiomatized, and there are only a few examples:
\begin{enumerate}
\item  The first two groups, considered in \cite{wa1}, are
$G=U_n$ and $G=O_n$. The symmetric group $G=S_n$ was considered
a few years later, in \cite{wa}. The corresponding quantum groups
have been intensively studied since then. However, the
``liberation'' aspect of the construction $G\to G^+$ was
understood only recently: the idea is that with $n\to\infty$ the
integral geometry of $G$ is governed by probabilistic
independence, while that of $G^+$, by freeness in the sense
of Voiculescu. This follows from \cite{bc1}, \cite{bc2},
\cite{cs}.

\item A number of subgroups $G\subset S_n$ have been considered
recently, most of them being symmetry groups of graphs, or of
other combinatorial structures. Some freeness appears here as
well, for instance in certain situations the operation
$G\to G^+$ transforms usual wreath products into free
wreath products. However, the main result so far goes
somehow in the opposite sense: for certain classes of groups, with
$n\to\infty$ we have $G^+=G$. See \cite{bb1}, \cite{bb2},
\cite{bbc}.
\end{enumerate}

Summarizing, the study of free quantum groups focuses on several
ad-hoc constructions $G\to G^+$. The following two questions
are of particular interest in connection with the axiomatization
problem: when does $G^+$ collapse to $G$? when is $G^+$ a ``true'' free version of $G$? There are already several
answers to these questions, most of them being of asymptotic
nature. Work here is in progress.

In this paper we study the free quantum group
associated to $G=S_4$. This corersponds to the Hopf algebra
$A_s(4)$ generated by $16$ abstract variables, subject to certain
relations, similar to those satisfied by the coordinates of
$S_4\subset U_4$. These relations, discovered in \cite{wa}, are
known under the name ``magic condition''.

We should mention first that for $G=S_n$ with $n=1,2,3$ we have
$A_s(n)=C(S_n)$, and there is no quantum group to be
studied. The other fact is that for $n\geq 5$ the algebra
$A_s(n)$ is too big, for instance it is not amenable, so a
detailed study here would be much more difficult, and would
require new ideas.

The interest in $A_s(4)$ comes from the fact that the space formed by $4$ points is the simplest
one to have a non-trivial quantum automorphism group. This suggests that $A_s(4)$ itself might be the
``simplest'' Hopf algebra, and we believe indeed that it is so. This algebra corresponds for instance to the situation ``index $4$, graph
$A_\infty$'', known from the work of Jones to be the basic one. See
\cite{b3}.

We have three results about $A_s(4)$:
\begin{enumerate}
\item An explicit realization of the magic condition is found in
\cite{bm}, by using the Pauli matrices. This gives a
representation $\pi:A_s(4)\to C(SU_2,M_4(\mathbb C))$,
shown there to be faithful in some weak sense. At that time there
was no technique for proving or disproving the fact that $\pi$ is
faithful. This problem can be solved now by using integration
techniques, inspired from \cite{bi}, \cite{co}, \cite{cs},
\cite{we} and introduced in \cite{bc1}, \cite{bc2}, and our result
is that $\pi$ is faithful. In other words, we have a model for the
abstract algebra $A_s(4)$. This is the second known model for a
free quantum algebra, the first one being a certain embedding
$A_u(2)\to C(\mathbb T)*C(SU_2)$, found in
\cite{b1}.

\item The model can be used for working out in detail the integral
geometry of the corresponding quantum group. As in our previous
work \cite{bc1}, \cite{bc2}, the problem that we consider is that of
computing laws of sums of diagonal coordinates $u_{ii}$. Each such
coordinate is a projection of trace $1/4$, and it is known from
\cite{b2} that the sum $u_{11}+u_{22}+u_{33}+u_{44}$ is free
Poisson on $[0,4]$. Our main result here is that the law of $u_{11}+u_{22}$ is an
average between a Dirac mass at $0$, and the Lebesgue measure on
$[0,2]$. In other words, modulo $0$ and rescalings, the averages
$M_s=(u_{11}+\ldots +u_{ss})/s$ with $s=1,2,4$ correspond to a Dirac
mass, a Lebesgue measure, and a free Poisson law.

\item We compute as well the laws of the variables interpolating between $M_1,M_2$ and $M_2,M_4$. The results here are quite technical, and the whole study belongs to a ``higher order'' problematics for free quantum groups. The idea is that one possible escape from asymptotic philosophy for $G=O_n$, probably adaptable to $G=U_n,S_n$ as well, would be the use of meander determinants of Di
Francesco \cite{df}. The other possible solution is via 
matrix models, and we use here the Pauli model in order for making
some advances.
\end{enumerate}

The paper is organized as follows: 1 is an introduction to $A_s(4)$, in 2-4 we discuss the Pauli representation, and in 5-9 we compute probability measures.

\section{Quantum permutation groups}

Let $A$ be a $C^*$-algebra. That is, we have a complex
algebra with a norm and an involution, such that the Cauchy sequences
converge, and $||aa^*||=||a||^2$.

The basic example is $B(H)$, the algebra of bounded operators on a
Hilbert space $H$. In fact, any $C^*$-algebra appears as
subalgebra of some $B(H)$.

The key example is $C(X)$, the algebra of continuous
functions on a compact space $X$. By a theorem of Gelfand, any
commutative $C^*$-algebra is of this form.

\begin{definition}
Let $A$ be a $C^*$-algebra.
\begin{enumerate}
\item A projection of $A$ is an element $p\in A$ satisfying $p^2=p=p^*$.
\item Two projections $p,q\in A$ are called orthogonal when
$pq=0$. \item A partition of the unity of $A$ is a finite set of projections of $A$, which are pairwise orthogonal, and which sum up to $1$.
\end{enumerate}
\end{definition}

A projection in $B(H)$ is an orthogonal projection $P(K)$, where
$K\subset H$ is a closed subspace. Orthogonality of projections
corresponds to orthogonality of subspaces, and partitions of unity
correspond to decompositions of $H$.

A projection in $C(X)$ is a characteristic function
$\chi(Y)$, where $Y\subset X$ is an open and closed subset.
Orthogonality of projections corresponds to disjointness of
subsets, and partitions of unity correspond to partitions of $X$.

\begin{definition}
A magic unitary is a square matrix $u\in M_n(A)$, whose rows
and columns are all partitions of unity in $A$.
\end{definition}

A magic unitary over $B(H)$ is of the form $P(K_{ij})$, with $K$
magic decomposition of $H$, in the sense that all rows and columns
of $K$ are decompositions of $H$. The basic examples here are of
the form $K_{ij}=\mathbb C\,\xi_{ij}$, where $\xi$ is a magic
basis of $H$, in the sense that all rows and columns of $\xi$ are
bases of $H$.

A magic unitary over $C(X)$ is of the form $\chi(Y_{ij})$,
with $Y$ magic partition of $X$, in the sense that all rows and
columns of $Y$ are partitions of $X$.  The key example here comes
from a finite group $G$ acting on a finite set $X$: the
characteristic functions $\chi_{ij}=\left\{\sigma\in G\mid
\sigma(j)=i\right\}$ form a magic unitary over $C(G)$.

In the particular case of the symmetric group $S_n$ acting on
$\{1,\ldots ,n\}$, we have the following presentation result,
which follows from the Gelfand theorem:

\begin{theorem}
$C(S_n)$ is the commutative $C^*$-algebra
generated by $n^2$ elements $\chi_{ij}$, with relations making
$(\chi_{ij})$ a magic unitary matrix. Moreover, the maps
\begin{eqnarray*}
\Delta(\chi_{ij})&=&\sum \chi_{ik}\otimes \chi_{kj}\cr
\varepsilon(\chi_{ij})&=&\delta_{ij}\cr S(\chi_{ij})&=&\chi_{ji}
\end{eqnarray*}
are the comultiplication, counit and antipode of $C(S_n)$.
\end{theorem}

In other words, when regarding $S_n$ as an algebraic group, the
relations satisfied by the $n^2$ coordinates are those expressing
magic unitarity.

We are interested in the algebra of ``free coordinates'' on $S_n$.
This is obtained by removing commutativity in the above
presentation result:

\begin{definition}
$A_s(n)$ is the $C^*$-algebra generated by $n^2$
elements $u_{ij}$, with relations making $(u_{ij})$ a magic
unitary matrix. The maps
\begin{eqnarray*}
\Delta(u_{ij})&=&\sum u_{ik}\otimes u_{kj}\cr
\varepsilon(u_{ij})&=&\delta_{ij}\cr S(u_{ij})&=&u_{ji}
\end{eqnarray*}
are called the comultiplication, counit and antipode of $A_s(n)$.
\end{definition}

This algebra, discovered by Wang in \cite{wa}, fits into
Woronowicz's quantum group formalism in \cite{wo}. In fact, the
quantum group $S_n^+$ defined by the heuristic formula $A_s(n)=C(S_n^+)$ is a free analogue of the symmetric group
$S_n$. This quantum group doesn't exist of course: the idea is
just that various properties of $A_s(n)$ can be expressed in terms
of it. As an example, the canonical quotient map $A_s(n)\to C(S_n)$
should be thought of as coming from an embedding $S_n\subset S_n^+$.

\begin{proposition}
For $n=1,2,3$ we have $A_s(n)=C(S_n)$.
\end{proposition}

This follows from the fact that for $n\leq 3$, the entries of a
$n\times n$ magic unitary matrix have to commute. Indeed, at $n=2$
the matrix must be of the form
$$u=\begin{pmatrix}p&1-p\cr 1-p&p\end{pmatrix}$$
with $p$ projection, and entries of this matrix commute. For $n=3$
see \cite{wa}.

This result is no longer true for $n=4$, where more complicated
examples of magic unitary matrices are available, for instance via
diagonal concatenation:
$$u=\begin{pmatrix}
p&1-p&0&0\cr 1-p&p&0&0\cr 0&0&q&1-q\cr 0&0&1-q&q
\end{pmatrix}$$

This example shows that the algebra generated by $p,q$ is a
quotient of $A_s(4)$. With $pq\neq qp$ we get that $A_s(4)$ is
non-commutative, hence bigger that $C(S_4)$.

We have the following result, proved in \cite{b2}, \cite{b3},
\cite{wa}.

\begin{theorem}
The algebra $A_s(4)$ has the following properties:
\begin{enumerate}
\item It is non-commutative, and infinite dimensional. \item It is
amenable in the discrete quantum group sense. \item The fusion
rules are the same as those for $SO_3$.
\end{enumerate}
\end{theorem}

We would like to end with a comment about the algebra $A_s(n)$
with $n\geq 5$, which won't appear in what follows. As in the case
$n=4$, this algebra is non-commutative, infinite dimensional, and
has same fusion rules as $SO_3$.

The subtlety comes from the fact that the irreducible
corepresentations have dimensions bigger than those of $SO_3$,
and this makes this algebra non-amenable. In fact, much worse is
expected to be true: there is evidence from \cite{vv} that the
reduced version should be simple, and that the corresponding von
Neumann algebra should be a prime ${\rm II}_1$-factor. In other
words, this algebra has bad analytic properties.

\section{The Pauli representation}

The purpose of this paper is to provide a detailed analytic
description of $A_s(4)$. We use an explicit matrix model, coming
from the Pauli matrices:
$$c_1=\begin{pmatrix}1&0\cr 0&1\end{pmatrix}\hskip 5mm
c_2=\begin{pmatrix}i&0\cr 0&-i\end{pmatrix}\hskip 5mm
c_3=\begin{pmatrix}0&1\cr -1&0\end{pmatrix}\hskip 5mm
c_4=\begin{pmatrix}0&i\cr i&0\end{pmatrix}$$

These are elements of $SU_2$. In fact, $SU_2$ consists of linear
combinations of Pauli matrices, with points on the real sphere
$S^3$ as coefficients:
$$SU(2)=\left\{\sum x_ic_i\,\big\vert \sum|x_i|^2=1\right\}$$

The Pauli matrices multiply according to the formulae for
quaternions:
$$c_2^2=c_3^2=c_4^2=-1$$
$$c_2c_3=-c_3c_2=c_4$$
$$c_3c_4=-c_4c_3=c_2$$
$$c_4c_2=-c_2c_4=c_3$$

The starting remark is that the Pauli matrices form an orthonormal
basis of $M_2(\mathbb C)$, with respect to the scalar product
$<a,b>=tr(b^*a)$. Moreover, the same is true if we multiply them
to the left or to the right by an element of $SU_2$.

We can formulate this fact in the following way.

\begin{theorem}
For any $x\in SU_2$ the elements $\xi_{ij}=c_ixc_j$ form a magic basis of $M_2({\mathbb C})$, with respect to the sclar product $<a,b>={\rm tr}(b^*a)$.
\end{theorem}

We fix a Hilbert space isomorphism $M_2({\mathbb C})\simeq
{\mathbb C}^4$, and we use the corresponding identification of
operator algebras $B(M_2({\mathbb C}))\simeq M_4({\mathbb C})$.

Associated to each element $x\in SU_2$ is the representation of $A_s(4)$ mapping $u_{ij}$ to the rank one projection on $c_ixc_j$:
$$\pi_x:A_s(4)\to M_4(\mathbb C)$$

This representation depends on $x$. For getting a faithful
representation, the idea is to regard all these representations as
fibers of a single representation.

\begin{definition}
The Pauli representation of $A_s(4)$ is the map
$$\pi:A_s(4)\to C(SU_2, M_4({\mathbb C}))$$
mapping $u_{ij}$ to the function $x\to$ rank one projection on $c_ixc_j$.
\end{definition}

As a first remark, in this statement $SU_2$ can be replaced by $PU_2=SO_3$. For reasons that will become clear later on, we
prefer to use $SU_2$.

This representation is introduced in \cite{bm}, with the main
result that it is faithful in some weak sense. In what follows we
prove that $\pi$ is faithful, in the usual sense.

Consider the natural linear form on the algebra on the right:
$$\int\varphi=\int_{SU_2} {\rm tr}\left(\varphi(x)\right)\, dx$$

We use the following analytic formulation of faithfulness.

\begin{proposition}
The representation $\pi$ is faithful provided that
$$\int u_{i_1j_1}\ldots u_{i_kj_k}=
\int \pi_{i_1j_1}\ldots \pi_{i_kj_k}$$ for any choice of $k$ and
of various $i,j$ indices, where $\pi_{ij}=\pi(u_{ij})$.
\end{proposition}

\begin{proof}
The condition in the statement is that
$$\int a=\int\pi(a)$$
for any product $a$ of generators $u_{ij}$. By linearity and
density this formula holds on the whole algebra $A_s(4)$. In other
words, the following diagram commutes:
$$\begin{matrix}
A_s(4)&\ &\displaystyle{\mathop{\longrightarrow}^ {\int}}&\
&{\mathbb C}\cr \ \cr \pi\downarrow&\ &\ &\ &\uparrow {\rm tr}\cr
\ \cr C\left(SU_2,M_4({\mathbb C})\right)&\ &
\displaystyle{\mathop{\longrightarrow}^ {\int}}&\ &M_4({\mathbb
C})
\end{matrix}$$

On the other hand, we know that $A_s(4)$ is amenable in the discrete quantum group sense. This means that its Haar functional
is faithful:
$$a\neq 0\Longrightarrow \int aa^*>0$$

Assume now that we have $\pi(a)=0$. This implies $\pi(aa^*)=0$,
and commutativity of the above diagram gives $\int aa^*=0$. Thus
$a=0$, and we are done.
\end{proof}

We compute now the integral on the right in the above statement.

Consider the canonical action of $SU_2$ on the algebra
$M_2(\mathbb C)^{\otimes k}$, obtained as $k$-th tensor power of
the adjoint action on $M_2(\mathbb C)$:
$$\alpha_x(a_1\otimes\ldots\otimes a_k)=xa_1x^*\otimes\ldots\otimes xa_kx^*$$

The following map will play an important role throughout this paper.

\begin{definition}
We define a linear map $R:M_2(\mathbb C)^{\otimes k}\to M_2(\mathbb C)^{\otimes k}$ by
$$R(c_{i_1}\otimes\ldots\otimes  c_{i_k})=\frac{1}{2}(c_{i_1}c_{i_2}^*
\otimes c_{i_2}c_{i_3}^*\otimes\ldots\otimes c_{i_k}c_{i_1}^*)$$
with the convention that for $k=1$ we have $R(c_i)=c_ic_i^*/2=1/2$.
\end{definition}

To any multi-index $i=(i_1,\ldots ,i_k)$ we associate the following element:
$$c_i=c_{i_1}\otimes\ldots\otimes c_{i_k}$$

With these notations, we have the following result.

\begin{proposition}
We have the formula
$$\int \pi_{i_1j_1}\ldots \pi_{i_kj_k}=
<R^*ER(c_j),c_i>$$ where $E$ is
the expectation under the canonical action of $SU_2$.
\end{proposition}

\begin{proof}
We have the following computation, where $P(\xi)$ denotes the rank
one projection onto a vector $\xi$:
\begin{eqnarray*}
\int \pi_{i_1j_1}\ldots \pi_{i_kj_k} &=& \int
\pi(u_{i_1j_1})\ldots \pi(u_{i_kj_k})\cr &=&\int
P(c_{i_1}xc_{j_1})\ldots P(c_{i_k}xc_{j_k})\cr &=&\int_{SU_2}{\rm
tr}\left( P(c_{i_1}xc_{j_1})\ldots P(c_{i_k}xc_{j_k})\right)\,dx
\end{eqnarray*}

We use now the following elementary formula, valid for any
sequence of norm one vectors $\xi_1,\ldots,\xi_k$ in a Hilbert
space:
\begin{eqnarray*}
{\rm Tr}\left(P(\xi_1)\ldots P(\xi_k)\right)
&=&<\xi_1,\xi_2><\xi_2,\xi_3>\ldots <\xi_k,\xi_1>
\end{eqnarray*}

In our situation these vectors are in fact matrices, and their
scalar products are given by $<\xi,\eta>={\rm tr}(\xi\eta^*)$.
This gives the following formula:
\begin{eqnarray*}
\int \pi_{i_1j_1}\ldots \pi_{i_kj_k} &=&
\frac{1}{4}\int_{SU_2}<c_{i_1}xc_{j_1},c_{i_2}xc_{j_2}>\ldots
<c_{i_k}xc_{j_k},c_{i_1}xc_{j_1}>\,dx\cr &=&
\frac{1}{4}\int_{SU_2}
 {\rm
tr}(c_{i_1}xc_{j_1}c_{j_2}^*x^*c_{i_2}^*)\ldots {\rm
tr}(c_{i_k}xc_{j_k}c_{j_1}^*x^*c_{i_1}^*)\,dx \cr
&=&\frac{1}{4}\int_{SU_2}
 {\rm
tr}(c_{i_2}^*c_{i_1}xc_{j_1}c_{j_2}^*x^*)\ldots {\rm
tr}(c_{i_1}^*c_{i_k}xc_{j_k}c_{j_1}^*x^*)\,dx
\end{eqnarray*}

We use now the formula $c_s^*=\pm c_s$, valid for all matrices
$c_s$. The minus signs can be rearranged, and the computation can
be continued as follows:
\begin{eqnarray*}
\int \pi_{i_1j_1}\ldots \pi_{i_kj_k} &=&\frac{1}{4}\int_{SU_2}
 {\rm
tr}(c_{i_2}c_{i_1}^*xc_{j_1}c_{j_2}^*x^*)\ldots {\rm
tr}(c_{i_1}c_{i_k}^*xc_{j_k}c_{j_1}^*x^*)\,dx\cr &=&
\frac{1}{4}\int_{SU_2}{\rm
tr}\left(c_{i_2}c_{i_1}^*xc_{j_1}c_{j_2}^*x^*\otimes\ldots\otimes
c_{i_1}c_{i_k}^*xc_{j_k}c_{j_1}^*x^*\right)\,dx\cr
&=&\int_{SU_2}{\rm tr}\left( R(c_i)^*\alpha_x(R(c_j))\right)\,dx
\end{eqnarray*}

We can interchange the trace and integral signs:
\begin{eqnarray*}
\int \pi_{i_1j_1}\ldots \pi_{i_kj_k}&=&{\rm tr} \left(\int_{SU_2}
R(c_i)^*\alpha_x(R(c_j))\, dx\right)\cr &=&{\rm tr} \left(
R(c_i)^*\int_{SU_2} \alpha_x(R(c_j))\,dx\right)
\end{eqnarray*}

Now acting by group elements, then integrating, is the same as
projecting onto fixed points, so the computation can be continued
as follows:
\begin{eqnarray*}
\int \pi_{i_1j_1}\ldots \pi_{i_kj_k}&=&{\rm
tr}\left(R(c_i)^*ER(c_j)\right)\cr &=&<ER(c_j),R(c_i)>\cr
&=&<R^*ER(c_j),c_i>
\end{eqnarray*}

This finishes the proof.
\end{proof}

\section{Some technical results}

We know from Proposition 2.1 that the faithfulness of the Pauli representation is equivalent to a certain equality of integrals. Moreover, one of these integrals can be computed by using the Weingarten formula in \cite{bc2}. As for the other integral, Proposition 2.2 shows that this can be computed as well, provided we have enough information about the operator $R$ from Definition 2.2.

In this section we work out a number of technical properties of this operator $R$. For this purpose, we need first understand some aspects of the structure of the algebra of fixed points under the diagonal adjoint action of
$SU_2$.

\begin{lemma}
$f=\sum_{i=1}^4c_i\otimes c_i^*$
is invariant under the action of $SU_2$.
\end{lemma}

\begin{proof}
We have $c_i=-c_i^*$ for $i=2,3,4$, so the element in the statement is:
$$f=2(c_1\otimes c_1)-\sum_{i=1}^4 c_i\otimes c_i$$

Since $c_1\otimes c_1$ is invariant under $SU_2$, what is left to prove is that the sum on the right is invariant as well. This sum, viewed as a matrix, is:
$$\sum_{i=1}^4c_i\otimes c_i=\begin{pmatrix}
 0&0&0&0\\
0&2&-2&0\\
0&-2&2&0\\
0&0&0&0
\end{pmatrix}$$

But this matrix is $4$ times the projection onto the determinant subspace of 
$\mathbb{C}^{2}\otimes\mathbb{C}^2$, which is invariant under the action of $SU_2$. This concludes the proof. 
\end{proof}

Let $NC(k)$ be the set of noncrossing partitions of $\{1,\ldots,k\}$. Given a partition $p\in NC(k)$ and a multi-index $j=(j_1,\ldots,j_k)$, we can plug $j$ into $p$ in the obvious way, and we define a number $\delta_{pj}\in\{0,1\}$ in the following way: $\delta_{pj}=1$ if any block of $p$ contains identical indices of $j$, and $\delta_{pj}=0$ if not. See \cite{bc2}. 

To any $p\in NC(k)$ we associate an element of $M_2(\mathbb C)^k$, in the following way:
$$c_p=\sum_j\delta_{pj}\,c_j$$

Our next goal is to prove that $R(c_p)$ is invariant under the action of $SU_2$. We discuss first the case of the trivial partition, $0_k=\{\{1\},\ldots ,\{k\}\}$.

\begin{lemma}
$R(c_{0_k})$ is invariant under the action of $SU_2$.
\end{lemma}

\begin{proof}
The partition $p=0_k$ has the particular property that $\delta_{pj}=1$ for any multi-index $j=(j_1,\ldots ,j_k)$. This gives the following formula:
\begin{eqnarray*}
R(c_{0_k})
&=&R\left(\sum_jc_j\right)\cr
&=&R\left(\sum_{j_1\ldots j_k}c_{j_1}\otimes\ldots\otimes c_{j_k}\right)\cr
&=&\frac{1}{2}\sum_{j_1\ldots j_k} c_{j_1}c^*_{j_2}\otimes
c_{j_2}c^*_{j_3}\otimes\ldots\otimes c_{j_k}c_{j_1}^*
\end{eqnarray*}

To any multi-index $j=(j_1,\ldots ,j_k)$ we associate a multi-index $i=(i_1,\ldots ,i_{k-1})$ in the following way: $i_1$ is such that $c_{i_1}=\pm c_{j_1}c_{j_2}^*$, $i_2$ is such that $c_{i_2}=\pm c_{j_1}c_{j_3}^*$, and so on up to the last index $i_{k-1}$, which is such that $c_{i_{k-1}}=\pm c_{j_1}c_{j_k}^*$. 

With this notation, we have the following formulae, where the possible dependences between the various  $\pm$ signs are not taken into account:
\begin{eqnarray*}
c_{i_1}&=&\pm c_{j_1}c_{j_2}^*\cr
c_{i_1}^*c_{i_2}&=&(\pm c_{j_1}c_{j_2}^*)^*(\pm c_{j_1}c_{j_3}^*)=\pm c_{j_2}c_{j_3}^*\cr
c_{i_2}^*c_{i_3}&=&(\pm c_{j_1}c_{j_3}^*)^*(\pm c_{j_1}c_{j_4}^*)=\pm c_{j_3}c_{j_4}^*\cr
\ldots&\ldots&\ldots\cr
c_{i_{k-2}}^*c_{i_{k-1}}&=&(\pm c_{j_1}c_{j_{k-1}}^*)^*(\pm c_{j_1}c_{j_k}^*)=\pm c_{j_{k-1}}c_{j_k}^*\cr
c_{i_{k-1}}^*&=&(\pm c_{j_1}c_{j_k}^*)^*=\pm c_{j_k}c_{j_1}^*
\end{eqnarray*}

By taking the tensor product of all these formulae, we get:
$$c_{i_1}\otimes c_{i_1}^*c_{i_2}\otimes \ldots\otimes c_{i_{k-2}}^*c_{i_{k-1}} \otimes c_{i_{k-1}}^*=\pm c_{j_1}c_{j_2}^*\otimes
c_{j_2}c^*_{j_3}\otimes \ldots \otimes c_{j_k}c_{j_1}^*$$

By applying the linear map given by $a_1\otimes\ldots \otimes a_k\to a_1\ldots a_k$ to both sides we see that the sign on the right is actually $+$. That is, we have:
$$c_{i_1}\otimes c_{i_1}^*c_{i_2}\otimes \ldots\otimes c_{i_{k-2}}^*c_{i_{k-1}} \otimes c_{i_{k-1}}^*=c_{j_1}c_{j_2}^*\otimes
c_{j_2}c^*_{j_3}\otimes \ldots \otimes c_{j_k}c_{j_1}^*$$

We recognize at right the basic summand in the formula of $R(c_{0_k})$. Now it follows from the definition of $i$ that summing the right terms over all multi-indices $j=(j_1,\ldots ,j_k)$ is the same as summing the left terms over all multi-indices $i=(i_1,\ldots,i_{k-1})$, then multiplying by $4$. Thus we get:
\begin{eqnarray*}
R(c_{0_k})
&=&\frac{4}{2}\sum_{i_1\ldots i_{k-1}}c_{i_1}\otimes c_{i_1}^*c_{i_2}\otimes \ldots \otimes c_{i_{k-2}}^*c_{i_{k-1}}\otimes c_{i_{k-1}}^*\cr
&=&2\sum_{i_1\ldots i_{k-1}}(c_{i_1}\otimes c_{i_1}^*)_{12}(c_{i_2}\otimes c_{i_2}^*)_{23}\ldots (c_{i_{k-1}}\otimes c_{i_{k-1}}^*)_{k-1,k}\cr
&=&2\left(\sum_{i_1}c_{i_1}\otimes c_{i_1}^*\right)_{12}\left(\sum_{i_2}c_{i_2}\otimes c_{i_2}^*\right)_{23}\ldots \left(\sum_{i_{k-1}}c_{i_{k-1}}\otimes c_{i_{k-1}}^*\right)_{k-1,k}\cr
&=&2f_{12}f_{23}\ldots f_{k-1,k}
\end{eqnarray*}

Now $f$ being invariant under $SU_2$, 
it is the same for each $f_{i,i+1}$. Since the invariants form an algebra, the above product is invariant. This concludes the proof.
\end{proof}

The aim of the next lemma is to prove that the previous lemma holds not only for
$c_0$, but also for $c_p$. 

We first introduce some notations. The Kreweras complement of a partition $p\in NC(k)$ is constructed as follows. Consider the ordered set $\{1,\ldots ,k\}$. At right of each index $i$ we put an index $i'$, as to get the following sequence of indices:
$$1,1',2,2',\ldots ,k,k'$$

The Kreweras complement $p^c$ is then the largest noncrossing partition of the new index set $\{1',\ldots ,k'\}$, such that the union of $p$ and $p^c$ is noncrossing.

\begin{definition}
For a noncrossing partition $p$ we use the notation
$$\omega (p)=2\,\prod_{i=1}^l\left(\prod_{p=1}^{k_i-1} f_{j_{i,p}j_{i,p+1}}\right)$$
where $\{j_{11}<\ldots <j_{1k_1}\},\ldots ,\{j_{l1}<\ldots <j_{lk_l}\}$ are the blocks of $p$, with the convention that for $k_i=1$ the product on the right is by definition $1$. 
\end{definition}

Here we use the element $f$ from Lemma 3.1, and the leg-numbering notation. The element $\omega(p)$ is well-defined, because the products on the right pairwise commute.

In particular in the case of the partition 
$1_k=\{\{1,\ldots ,k\}\}$, we have:
$$\omega(1_k)=2\prod_{i=1}^{k-1}f_{i,i+1}$$

Since each $f_{ij}$ is invariant under the adjoint action of $SU_2$ and
since the invariants form an algebra, $\omega (p)$ is also invariant under the adjoint action of $SU_2$.

\begin{lemma}
$R(c_p)$ is invariant under the action of $SU_2$, for any $p\in NC(k)$.
\end{lemma}

\begin{proof}
The idea is to generalize the proof of the previous lemma. In particular we need to generalize the key formula there, namely:
$$R(c_{0_k})=2f_{12}f_{23}\ldots f_{k-1,k}$$

We claim that for any noncrossing partition $p$, we have:
$$R(c_p)=\omega(p^c)$$

As a first verification, for the trivial partition $p=0_k$ we have $p^c=1_k$, and the formula $R(c_p)=\omega(p^c)$ follows from the above identities. Also, for the rough partition $p=1_k$ we have $p^c=0_k$, and the claimed equality follows from:
\begin{eqnarray*}
R(c_{1_k})
&=&\sum_i R(c_i\otimes\ldots\otimes c_i)\\
&=&\frac{1}{2}\sum_ic_ic_i^*\otimes\ldots\otimes c_ic_i^*\\
&=&\frac{1}{2}\,4(1\otimes\ldots\otimes 1)\\
&=&\omega(0_k)
\end{eqnarray*}

In the general case, we can use the same method as for Lemma 3.2. We have:
\begin{eqnarray*}
R(c_{p})
&=&R\left(\sum_j\delta_{pj}c_j\right)\cr
&=&R\left(\sum_{j_1\ldots j_k}
\delta_{pj}c_{j_1}\otimes\ldots\otimes c_{j_k}\right)\cr
&=&\frac{1}{2}\sum_{j_1\ldots j_k} \delta_{pj}c_{j_1}c^*_{j_2}\otimes
c_{j_2}c^*_{j_3}\otimes\ldots\otimes c_{j_k}c_{j_1}^*
\end{eqnarray*}

To any multi-index $j=(j_1,\ldots ,j_k)$ we associate 
a multi-index $i=(i_1,\ldots ,i_k)$ in the following way: $i_1$ is such that $c_{i_1}=\pm c_{j_1}c_{j_2}^*$, $i_2$ is such that $c_{i_2}=\pm c_{j_1}c_{j_3}^*$, and so on up to the last index $i_k$, which is such that $c_{i_{k}}=\pm c_{j_k}c_{j_1}^*$. 

With this notation, together with the observation that the product of all the above $\pm$ signs is actually $+$ (this can be checked as in the proof of Lemma 3.2), the above formula becomes:
$$R(c_p)=\frac{1}{2}\sum_j\delta_{pj}c_{i_1}\otimes\ldots\otimes c_{i_k}$$

Let $l$ be an index in $\{1,\,\ldots ,k\}$ and assume that 
it is a last element 
of a block of $p^c$.
Let $l_1<\ldots<l_x=l$ be the ordered enumeration of the elements 
of this block. 
We observe that: 
$$c_{j_{l_1}}\ldots c_{j_{l_x}}=1$$

For example,
if $l$ and $l+1$ are in the same block of $p$, then
 $\{l\}$ is a one-element block in $p^c$ and $i_l$ will be such
that $c_{i_l}=\pm c_{j_l}c_{j_{l+1}}^*=c_{j_l}c_{j_l}^*=1$.
The general case works by following the argument of Lemma 3.2.
We get:
$$\sum_i\delta_{pi}c_{j_{l_1}}\otimes \ldots 
\otimes c_{j_{l_x}}=R(c_{0_x})$$

The point is that $R(c_p)$ is the product of the above expressions, over the blocks of $p^c$. If we denote these blocks by $\{l_{11}<\ldots<l_{1x_1}\}$, \ldots , $\{l_{r1}<\ldots<l_{rx_r}\}$, we get, by using the leg-numbering notation:
\begin{eqnarray*}
R(c_p)
&=&2^{1-b}\sum_j\prod_{b=1}^rR(c_{0_{x_b}})_{l_{b1}\ldots l_{bx_b}}\\
&=&2^{1-b}\sum_j\prod_{b=1}^r\omega(1_{x_b})_{l_{b1}\ldots l_{bx_b}}\\
&=&\omega(p^c)
\end{eqnarray*}

To illustrate the above proof we perform explicitly the computation 
in the case of the following partition:
$$p=\{\{1,5\},\{2\},\{3,4\},\{6\}\}$$ 

The Kreweras complement of this partition is $p^c=\{\{1,2,4\},\{3\},\{5,6\}\}$, and the above method gives:
\begin{eqnarray*}
R(c_p)
&=&\frac{1}{2}\sum_{j_1j_2j_3j_4}c_{j_1}c_{j_2^*}\otimes c_{j_2}c_{j_3}^*\otimes 1\otimes c_{j_3}c_{j_1}^*\otimes c_{j_1}c_{j_4}^*\otimes c_{j_4}c_{j_1}^*\\
&=&\frac{1}{2}\,4\sum_{i_1i_2i_3}c_{i_1}\otimes c_{i_1}^*c_{i_2}\otimes 1\otimes c_{i_2}^*\otimes c_{i_3}\otimes c_{i_3}^*\\
&=&2\sum_{i_1i_2i_3}(c_{i_1}\otimes c_{i_1}^*c_{i_2}\otimes 1\otimes c_{i_2}^*\otimes 1\otimes 1)(1\otimes 1\otimes 1\otimes 1\otimes c_{i_3}\otimes c_{i_3}^*)\\
&=&2(f_{12}f_{24})f_{56}\\
&=&\omega(p^c)
\end{eqnarray*}

Now back to the general case, the formula $R(c_p)=\omega(p^c)$ shows that $R(c_p)$ is a product of certain elements $f_{ij}$ obtained from $f$ by acting on the various legs of $M_2(\mathbb C)^{\otimes k}$, and we can conclude as in the proof of the previous lemma.
\end{proof}

\section{Proof of faithfulness}

We are now in position of proving the main result. With the technical results from previous section, we can describe the operator $R^*ER$ from Proposition 2.1.

\begin{proposition}
We have $R^*ER(c_p)=c_p$.
\end{proposition}

\begin{proof}
The previous lemma gives the following formula:
\begin{eqnarray*}
<R^*ER(c_p),c_j>
&=&<R^*R(c_p),c_j>\cr
&=&<R(c_p),R(c_j)>\cr
&=&\sum_i\delta_{pi}<R(c_i),R(c_j)>
\end{eqnarray*}

On the other hand, we have from definitions:
$$<c_p,c_j>=\delta_{pj}$$

Since the elements $c_j$ span the ambient space, what is left to prove is:
$$\sum_i\delta_{pi}<R(c_i),R(c_j)>=\delta_{pj}$$

But this can be checked by direct computation. First, from the definition of $R$, we have the following formula:
\begin{eqnarray*}
<R(c_i),R(c_j)>
&=&\frac{1}{4}\,<c_{i_1}c_{i_2}^*\otimes\ldots\otimes
c_{i_k}c_{i_1}^*,c_{j_1}c_{j_2}^*\otimes\ldots\otimes
c_{j_k}c_{j_1}^*>\cr&=&\frac{1}{4}\,<c_{i_1}c_{i_2}^*,c_{j_1}c_{j_2}^*>
\ldots <c_{i_k}c_{i_1}^*,c_{j_k}c_{j_1}^*>
\end{eqnarray*}

In this formula all scalar products are $0$, $1$ and $-1$. Now assume that $i_1,\ldots ,i_k$ and $j_1$ are fixed. If the first scalar product is $\pm 1$ then $j_2$ is uniquely determined, then if the second scalar product is $\pm 1$ then $j_3$ is uniquely determined as well, and so on. Thus for all scalar products to be $\pm 1$, the multi-index $j$ is uniquely determined by the multi-index $i$, up to a
possible choice of the first index $j_1$.

Moreover, each choice of $j_1$ leads to a multi-index $j=(j_1,\ldots,j_k)$, such that all the scalar products are $\pm 1$. Indeed, once $j_2,\ldots,j_k$ are chosen as to satisfy $c_{i_r}c_{i_{r+1}}^*=\pm c_{j_r}c_{j_{r+1}}^*$ for $r=1,\ldots,k-1$, by multiplying all these formulae we get $c_{i_1}c_{i_k}^*=\pm c_{j_1}c_{j_k}^*$, which shows that the last scalar product is $\pm 1$ as well.

Summarizing, given a multi-index $i=(i_1,\ldots ,i_k)$ and a number $s\in\{1,2,3,4\}$, there is a unique multi-index $j=(j_1,\ldots,j_k)$ with $j_1=s$, such that all the above scalar products are $\pm 1$. We use the notation $j=i\oplus s$.

In the situation $j=i\oplus s$ we have $<R(c_i),R(c_j)>=\pm 1$, and by applying the linear map given by $a_1\otimes\ldots \otimes a_k\to a_1\ldots a_k$ to the formula $c_{i_1}c_{i_2}^*\otimes\ldots\otimes c_{i_k}c_{i_1}^*=\pm c_{j_1}c_{j_2}^*\otimes\ldots\otimes c_{j_k}c_{j_1}^*$ we see that the sign is actually $+$. Thus we have:
$$<R(c_i),R(c_j)>
=\begin{cases}
1/4\mbox{ if $j=i\oplus s$ for some $s\in\{1,2,3,4\}$}\\
0\mbox{ if not}
\end{cases}$$

We can come back now to the missing formula. We have:
$$\sum_i\delta_{pi}<R(c_i),R(c_j)>
=\frac{1}{4}\sum_{s=1}^4\delta_{pi}$$

We claim that we have $\delta_{pi}=\delta_{pj}$, for any partition $p$. Indeed,this follows from the fact that for $r<s$ we have $i_r=i_s$ if and only if the product of $c_{i_t}c_{i_{t+1}}^*$ over $t=r,\ldots,s-1$ equals $\pm 1$, and a similar statement holds for the multi-index $j$.

We can now conclude the proof. By using $\delta_{pi}=\delta_{pj}$, we get:
$$\sum_i\delta_{pi}<R(c_i),R(c_j)>
=\frac{1}{4}\sum_{s=1}^4\delta_{pj}=\delta_{pj}$$

But this is the formula that we wanted to prove, so we are done.
\end{proof}

\begin{theorem}
The Pauli representation of $A_s(4)$ is faithful.
\end{theorem}

\begin{proof}
We denote as usual by $c_1,\ldots,c_4$ the Pauli matrices, and we let $e_1,\ldots,e_4$ be the standard basis of $\mathbb C^4$. For a multi-index $i=(i_1,\ldots ,i_k)$, we set:
\begin{eqnarray*}
e_i&=&e_{i_1}\otimes\ldots\otimes e_{i_k}\cr
c_i&=&c_{i_1}\otimes\ldots\otimes c_{i_k}
\end{eqnarray*}

Each partition $p\in NC(k)$ creates two tensors, in the following
way:
\begin{eqnarray*}
e_p&=&\sum_i\delta_{pi}\,e_i\cr
c_p&=&\sum_i\delta_{pi}\,c_i\end{eqnarray*}

Consider the following $4^k\times 4^k$ matrices, with entries labeled by multi-indices $i,j$:
$$P_{ij}=tr\left(\int \pi_{i_1j_1}\ldots \pi_{i_kj_k}\right)$$
$$U_{ij}=\left(\int u_{i_1j_1}\ldots u_{i_kj_k}\right)$$

According to Proposition 2.2, we have the following formula:
$$P_{ij}=<R^*ER(c_j),c_i>$$

Let $\Phi$ be the linear map $(\mathbb{C}^4)^{\otimes k}\to M_2(\mathbb{C})^{\otimes k}$
given by $\Phi (e_i)=c_i$. The above equation can be rephrased as:
$$P_{ij}=<\Phi^*R^*ER\Phi(e_j),e_i>$$

According to Proposition 2.1, it is enough to prove that we have $P=U$. But $U$ is the orthogonal projection of $(\mathbb{C}^4)^{\otimes k}$ onto the following space (see \cite{bc2}):
$$S_e={\rm span}\{ e_p\mid p\in NC(k)\}$$

Thus what we have to prove is that $\Phi^*R^*ER\Phi$ is the projection onto $S_e$. But this is equivalent to proving that $R^*ER$ is the projection onto the following space: 
$$S_c={\rm span}\{ c_p\mid p\in NC(k)\}$$

We know from Proposition 4.1 that $R^*ER(c_p)=c_p$, for all $p\in NC(k)$. This implies that
$R^*ER$ restricted to $S_c$ is the identity.
This vector space has dimension the Catalan number $C_k$, and this is 
exactly the rank of the operator $E$. Therefore 
$S_c$ has to be the image of
$R^*ER$. Now since $R^*ER$ is self-adjoint, and is the identity on its image, it is the orthogonal projection onto $S_c$, and this concludes the proof.
\end{proof}

\section{Diagonal coefficients}

We can use the Pauli representation for computing laws of certain
diagonal coefficients of $u$, the fundamental corepresentation of
$A_s(4)$. These diagonal coefficients are introduced in
\cite{bc1}, \cite{bc2}. We are particularly interested in the
following elements:

\begin{definition}
Associated to $s=1,2,3,4$ is the average
$$M_s=\frac{1}{s}\,\left(u_{11}+u_{22}+\ldots +u_{ss}\right)$$
where $u$ is the fundamental corepresentation of $A_s(4)$.
\end{definition}

We use notions from free probability from \cite{vo}, \cite{vdn}.
Recall first that the free Poisson law is the following
probability measure on $[0,4]$:
$$\nu=\frac{1}{2\pi}\,\sqrt{4x^{-1}-1}\,dx$$

In this paper we use probability measures supported on $[0,1]$.
So, consider the corresponding compression of the free Poisson
law:
$$\nu_1=\frac{2}{\pi}\,\sqrt{x^{-1}-1}\,dx$$

We denote by $\lambda_1$ the Lebesgue measure on $[0,1]$.

\begin{lemma}
We have the following formulae, where $(a,b,c,d)\in S^3$:
\begin{enumerate}
\item ${\rm law}(a^2)=\nu_1$. \item ${\rm
law}(a^2+b^2)=\lambda_1$. \item ${\rm
law}(a^2+b^2+c^2+d^2)=\delta_1$.
\end{enumerate}
\end{lemma}

\begin{proof}
These results are all well-known:

(1) The variable $2a$ is semicircular on $[-2,2]$, as one can see
geometrically on $S^3$, or by using representation theory of
$SU_2$, so its square $4a^2$ is free Poisson.

(2) This follows from the fact that when projecting $S^3$ on the
unit disk, the uniform measure on $S^3$ becomes the uniform
measure on the unit disk.

(3) This follows from $a^2+b^2+c^2+d^2=1$.
\end{proof}

With these formulae in hand, we can compute the laws of
$M_1,M_2,M_4$. Modulo a Dirac mass at $0$, these are the Dirac,
Lebesgue and free Poisson laws:

\begin{theorem}
For $s=1,2,4$ we have the formula
$${\rm law}(M_s)=\left(1-\frac{s}{4}\right)\,\delta_0+
\frac{s}{4}\,\mu_s$$ where $\mu_1=\delta_1$, $\mu_2=\lambda_1$,
$\mu_4=\nu_1$.
\end{theorem}

\begin{proof}
The $s=1$ assertion is clear, and the $s=4$ one is known from
\cite{b2}. For reasons of uniformity of the proof, we will prove
these assertions as well.

Consider an element $x\in SU_2$, and write it as $x=x_1c_1+\ldots
+x_4c_4$. With the notations $c_i^2=\varepsilon_i$ and $c_ic_j
=\varepsilon_{ij}c_jc_i$, with
$\varepsilon_i,\varepsilon_{ij}\in\{\pm 1\}$, we have:
$$c_ixc_i=\sum_j \varepsilon_i\varepsilon_{ij}x_jc_j$$

This gives the following formula for the orthogonal projection
onto the vector $c_ixc_i$, denoted $\pi_{ii}$, as in previous
section:
$$<\pi_{ii}c_j,c_k>=\varepsilon_{ij}\varepsilon_{ik}x_j x_k$$

Now by using the more convenient notation $x=ac_1+bc_2+cc_3+dc_4$,
we get from the multiplication table of Pauli matrices the
following formulae:
$$\pi_{11}=\begin{pmatrix}
a^2&ab&ac&ad\cr
ab&b^2&bc&bd\cr
ac&bc&c^2&cd\cr
ad&bd&cd&d^2
\end{pmatrix}\hskip 12mm\pi_{22}=\begin{pmatrix} a^2&ab&-ac&-ad\cr ab&b^2&-bc&-bd\cr
-ac&-bc&c^2&cd\cr -ad&-bd&cd&d^2
\end{pmatrix}$$
$$\pi_{33}=\begin{pmatrix}
a^2&-ab&ac&-ad\cr
-ab&b^2&-bc&bd\cr
ac&-bc&c^2&-cd\cr
-ad&bd&-cd&d^2
\end{pmatrix}\hskip 8mm\pi_{44}=\begin{pmatrix} a^2&-ab&-ac&ad\cr -ab&b^2&bc&-bd\cr
-ac&bc&c^2&-cd\cr ad&-bd&-cd&d^2
\end{pmatrix}$$

We have $M_1=\pi_{11}$, and by making averages we get $M_2,M_4$:
$$M_2=\begin{pmatrix} a^2&ab&0&0\cr ab&b^2&0&0\cr 0&0&c^2&cd\cr
0&0&cd&d^2
\end{pmatrix}\hskip 10mm
M_4=\begin{pmatrix} a^2&0&0&0\cr 0&b^2&0&0\cr 0&0&c^2&0\cr
0&0&0&d^2
\end{pmatrix}$$

With these notations, we have to compute the following numbers:
$$\int {\rm tr}(M_s^k)=\int \left(\frac{u_{11}+\ldots +u_{ss}}{s}\right)^k
$$

We first compute the characteristic polynomial of each $M_s$:
\begin{eqnarray*}
\det(y-M_1)&=&y^3(y-1)\cr
\det(y-M_2)&=&y^2(y-a^2-b^2)(y-c^2-d^2)\cr
\det(y-M_4)&=&(y-a^2)(y-b^2)(y-c^2)(y-d^2)
\end{eqnarray*}

Thus we have the following diagonalisations:
\begin{eqnarray*}
M_1&\sim&{\rm diag}(0,0,0,1)\cr M_2&\sim&{\rm
diag}(0,0,a^2+b^2,c^2+d^2)\cr M_4&\sim&{\rm diag}(a^2,b^2,c^2,d^2)
\end{eqnarray*}

We take powers, we apply the trace, and we integrate:
\begin{eqnarray*}
\int {\rm tr} (M_1^k)&=&\frac{1}{4}\cr \int {\rm tr}
(M_2^k)&=&\frac{1}{4}\int(a^2+b^2)^k+(c^2+d^2)^k\cr \int {\rm tr}
(M_4^k)&=&\frac{1}{4}\int a^{2k}+b^{2k}+c^{2k}+d^{2k}
\end{eqnarray*}

By symmetry reasons we have:
\begin{eqnarray*}
\int {\rm tr} (M_1^k)&=&\frac{1}{4}\cr \int {\rm tr}
(M_2^k)&=&\frac{1}{2}\int(a^2+b^2)^k\cr \int {\rm tr}
(M_4^k)&=&\int a^{2k}
\end{eqnarray*}

The result follows now from Lemma 5.1.
\end{proof}

We would like to end this section with a technical result, to be used later on.

\begin{proposition}
The Cauchy transforms for $M_1,M_2,M_4$ are:
\begin{eqnarray*}
G_1(\xi)&=&\frac{1}{\xi}+\frac{1}{4(\xi^2-\xi)}\cr
G_2(\xi)&=&\frac{1}{2}\left(\frac{1}{\xi}-\log\left(1-\frac{1}{\xi}\right)\right)\cr
G_4(\xi)&=&2\left(1-\sqrt{1-\frac{1}{\xi}}\right)
\end{eqnarray*}
\end{proposition}

\begin{proof}
This follows from the above formulae.
\end{proof}

\section{Numeric results}

In previous section we computed the law of the average $M_s$, with $s=1,2,4$. The same proof doesn't apply to the missing variable $M_3$, because the corresponding matrix cannot be diagonalized explicitly. This technical problem is to be related to the general principle ``the sphere cannot be cut in three parts''.

The undiagonalizable (or cutting) matrix is as follows:

\begin{proposition}
$M_3$ has the same law as the random matrix
$$M_3=\frac{1}{3}\begin{pmatrix}
3a^2&-ab&-ac&-ad\cr -ab&3b^2&-bc&-bd\cr -ac&-bc&3c^2&-cd\cr
-ad&-bd&-cd&3d^2
\end{pmatrix}$$
depending on $x=(a,b,c,d)$ on the sphere $S^3$.
\end{proposition}

\begin{proof}
By symmetry reasons $M_3$ has the same law as
$(u_{22}+u_{33}+u_{44})/3$, which has in turn the same law as the
following matrix:
$$M_3=\frac{1}{3}\left(\pi_{22}+\pi_{33}+\pi_{44}\right)$$

By using the formulae of $\pi_{ii}$, we get the matrix in the
statement.
\end{proof}

Observe that $M_3$ has trace $1$, and has
$(a^{-1},b^{-1},c^{-1},d^{-1})$ as $0$-eigenvector. Thus the
characteristic polynomial of $M_3$ is of the following form:
$$Q(y)=y^4-y^3+Ky^2-Ly$$

An explicit computation gives the following formulae for $K,L$:
$$K=\frac{8}{9}\left(a^2b^2+a^2c^2+a^2d^2+
b^2c^2+b^2d^2+c^2d^2\right)$$
$$L=\frac{16}{27}\left(a^2b^2c^2+a^2b^2d^2+
a^2c^2d^2+b^2c^2d^2\right)$$

In principle, this can be used for computing the Cauchy transform:
$$G(\xi)=\frac{1}{4}\int\frac{Q'(\xi)}{Q(\xi)}\,dx$$

The corresponding integration problem on $S^3$ looks
particularly difficult, and we don't know how to solve it. However, we did a lot of related abstract or numeric computations, and we have the following result:

\begin{theorem}
The first moments of $N_3=3M_3$ are given by:

\begin{tabular}
[t]{|l|l|l|l|l|l|l|l|}\hline 
{\rm order}&$1$&$2$&$3$&$4$&$5$&$6$&$\dots$\\ \hline 
{\rm moment}&$3/4$&$5/4$&$5/2$&$109/20$&$25/2$&$4157/140$&$\dots$\\
\hline
\end{tabular}
\vskip2mm

\hfill\begin{tabular}
[t]{|l|l|l|l|l|l|l|l|l|l|}\hline 
$\dots$&$7$&$8$&$9$\\ \hline 
$\dots$&$1449/20$&$75877/420$&$64223/140$\\
\hline
\end{tabular}
\vskip2mm 
\end{theorem}

\begin{proof}
We have to compute the moments of the matrix $N_3(a,b,c,d)$,
where $(a,b,c,d)$ is uniformly distributed along the sphere $S^3$.

For this, let $(a',b',c',d')$ be independent standard 
Gaussian variables. We have the following equality of joint distributions:
$${\rm law}(a',b',c',d')={\rm law}(\rho a,\rho b,\rho c,\rho d)$$
where $\rho$ is the positive square root of a standard chi-square distribution of parameter $4$, independent from $(a,b,c,d)$.

Consider now the random matrix $N_3(a',b',c',d')$, obtained by replacing $(a,b,c,d)$ with $(a',b',c',d')$. Since $\rho$ is a real variable independent from $(a,b,c,d)$, whose even moments are the numbers $2^k(k+1)!$, we have the following computation:
\begin{eqnarray*}
\int N_3(a',b',c',d')^k&=&\int N_3(\rho a,\rho b,\rho c,\rho d)^k\cr
&=&\int\rho^{2k}N_3(a,b,c,d)^k\cr
&=&\int\rho^{2k}\int N_3(a,b,c,d)^k\cr
&=&2^k(k+1)!\int N_3(a,b,c,d)^k
\end{eqnarray*}

The integrals on the left can be computed with Maple, and we get the result.
\end{proof}

The above theorem shows that our matrix model has also a computational interest: indeed, a direct attempt to make the above computations only with the 
Weingarten function summation formulae (or even with the matrix model but without the above Gaussianization trick) could not yield more than 5-6 moments whereas the method described in the proof of the above theorem could easily yield up to 15 moments.

Yet, these computations don't shed any light on where these moments come from: for instance the number 64223, appearing in the above table, is prime.

As a last comment, we have the following negative result.

\begin{proposition}
${\rm law}(M_3)\neq \frac{1}{4}\,\delta_0+\frac{3}{4}\,{\rm law}\left(a^2+b^2+c^2\right)$.
\end{proposition}

Indeed, the second moments of these laws are respectively
$5/36$ and $15/32$. This can be checked by a routine computation,
and contradicts what one might want to conjecture, after a quick
comparison of Lemma 5.1 and Theorem 5.1.

\section{Lebesgue-Dirac interpolation}

In this section and in the next one we perform some technical computations. Our motivation is as follows. Consider the variables $u_{11}+\ldots +u_{ss}$ with $s=1,2,4$, whose laws are known. The variable $u_{11}+u_{22}+u_{33}$ has the same law as:
$$u_{22}+u_{33}+u_{44}=(u_{11}+u_{22}+u_{33}+u_{44})-u_{11}$$

Thus, we are in front of the following problem: we know how to compute the laws of 3 variables, and we want to compute the law of a certain 4-th variable, belonging the same projective plane. This can be regarded as being part of the more general problem of finding the law of an arbitrary point in the plane. 

We will work out here two simple computations in this sense.

We know how to compute the laws of the following two elements:
\begin{eqnarray*}
w_0&=&(u_{11}+u_{22})/2\cr
w_1&=& u_{11}
\end{eqnarray*}

We can consider the following element, interpolating between them:

\begin{definition}
To any real number $t$ we associate the element
$$w_t=\frac{1+t}{2}\cdot u_{11}+\frac{1-t}{2}\cdot u_{22}$$
where $u$ is the fundamental corepresentation of $A_s(4)$.
\end{definition}

By using the matrix formulae in section 5, the characteristic polynomial of $w_t$ is given by $P(y)=y^2Q(y)$, where $Q$ is the following degree two polynomial:
$$Q(y)=y^2-y+(1-t^2)(a^2+b^2)(c^2+d^2)$$

In other words, $w_t$ has a double $0$ eigenvalue, and a $2\times
2$ matrix block. The law of the matrix block can be computed by using the following lemma.

\begin{lemma}
The Cauchy transform of a $2\times 2$ random matrix $M$ having
characteristic polynomial $Q(y)=y^2-By+C$ is given by the
following formula:
$$G(\xi)=\frac{1}{\xi}+\frac{1}{2\xi}\sum_{p+q>0}\frac{(-1)^p}
{\xi^{2p+q}}\cdot\frac{2p+q}{p+q}\begin{pmatrix}p+q\cr q
\end{pmatrix}
\int B^{q}C^p$$
\end{lemma}

\begin{proof}
The eigenvalues of $M$ are the roots of $Q$, so we get:
\begin{eqnarray*}
\int {\rm tr}(M^k)
&=&2^{-k-1}\int (B+\sqrt{B^2-4C})^k+ (B-\sqrt{B^2-4C})^k\cr
&=&2^{-k}\sum_{n=0}^{k/2}
\begin{pmatrix}k\cr 2n\end{pmatrix}\int B^{k-2n}(B^2-4C)^n\cr
&=&2^{-k}\sum_{n=0}^{k/2}\sum_{p=0}^n (-4)^p\begin{pmatrix}k\cr
2n\end{pmatrix}\begin{pmatrix} n\cr p\end{pmatrix}\int B^{k-2p}C^p
\end{eqnarray*}

Here a sum from $0$ to a real number $r$ means by definition sum
from $0$ to the integral part of $r$. We compute now the Cauchy
transform:
\begin{eqnarray*}
G(\xi)&=&\sum_{k=0}^\infty \xi^{-k-1}\int{\rm tr}(M^k)\cr
&=&\frac{1}{\xi}\sum_{k=0}^\infty\sum_{n=0}^{k/2}\sum_{p=0}^n
\frac{(-4)^p}{(2\xi)^{k}}\begin{pmatrix}k\cr
2n\end{pmatrix}\begin{pmatrix} n\cr p\end{pmatrix}\int B^{k-2p}C^p
\end{eqnarray*}

We make the replacements $n=p+m$ and $k=2p+q$:
\begin{eqnarray*}
G(\xi)&=& \frac{1}{\xi}\sum_{p,q=0}^\infty\sum_{m=0}^{q/2}
\frac{(-4)^p}{(2\xi)^{2p+q}}\begin{pmatrix}2p+q\cr
2p+2m\end{pmatrix}\begin{pmatrix} p+m\cr p\end{pmatrix}\int
B^{q}C^p
\end{eqnarray*}

We use now the following standard identity, valid for $p+q>0$:
$$\sum_{m=0}^{q/2}\begin{pmatrix}2p+q\cr
2p+2m\end{pmatrix}\begin{pmatrix} p+m\cr p\end{pmatrix}
=2^{q-1}\,\frac{2p+q}{p+q}\begin{pmatrix}p+q\cr q
\end{pmatrix}$$

This gives the following formula:
\begin{eqnarray*}
G(\xi)&=&\frac{1}{\xi}+\frac{1}{\xi}\sum_{p+q>0} \frac{(-4)^p\cdot
2^{q-1}}{(2\xi)^{2p+q}}\cdot\frac{2p+q}{p+q}\begin{pmatrix}p+q\cr
q
\end{pmatrix}
\int B^{q}C^p\cr
&=&\frac{1}{\xi}+\frac{1}{2\xi}\sum_{p+q>0}\frac{(-1)^p}
{\xi^{2p+q}}\cdot\frac{2p+q}{p+q}\begin{pmatrix}p+q\cr q
\end{pmatrix}
\int B^{q}C^p
\end{eqnarray*}

This gives the formula in the statement.
\end{proof}

\begin{theorem}
The Cauchy transform of the law of $w_t$ is given by
$$G(\xi)=\frac{1}{2\xi}+\frac{1-2\xi}{\xi-\xi^2}\cdot\frac{{\rm arcsinh}(\sqrt{x}/2)}{\sqrt{4x+x^2}}$$
where $x=(1-t^2)/(\xi^2-\xi)$. 
\end{theorem}

\begin{proof}
We recall that the characteristic polynomial of $w_t$ is given by the formula
$P(y)=y^2Q(y)$, where $Q$ is the following degree two polynomial:
$$Q(y)=y^2-y+(1-t^2)(a^2+b^2)(c^2+d^2)$$

Thus the Cauchy transform of $w_t$ is the average between $1/\xi$ and
the Cauchy transform $G_B$ of the corresponding $2\times 2$ matrix block:
$$G(\xi)=\frac{1}{2}\left(\frac{1}{\xi}+G_B(\xi)\right)$$

We apply Lemma 7.1 with the above characteristic polynomial, which
is of the form $Q(y)=y^2-y+sD$, with $s=1-t^2$ and
$D=(a^2+b^2)(c^2+d^2)$. We get:
$$G_B(\xi)=\frac{1}{\xi}+\frac{1}{2\xi}
\sum_{p+q>0}\frac{(-s)^p}
{\xi^{2p+q}}\cdot\frac{2p+q}{p+q}\begin{pmatrix}p+q\cr q
\end{pmatrix}
\int D^p$$

Thus we get the following formula:
$$G(\xi)=\frac{1}{\xi}+\frac{1}{4\xi}\sum_{p+q>0}
\frac{(-s)^p}{\xi^{2p+q}}
 \cdot\frac{2p+q}{p+q}\begin{pmatrix}p+q\cr q
\end{pmatrix}
\int D^p$$

We can get rid of $p+q>0$ by using the Cauchy transform $G_0$:
\begin{eqnarray*}
G(\xi)&=&G_0(\xi)+\frac{1}{4\xi}\sum_{p=1}^\infty
\sum_{q=0}^\infty\frac{(-s)^p}{\xi^{2p+q}}
 \cdot\frac{2p+q}{p+q}\begin{pmatrix}p+q\cr q
\end{pmatrix}
\int D^p\cr &=&G_0(\xi)+\frac{1}{4\xi}\sum_{p=1}^\infty
\left(-\frac{s}{\xi}\right)^p\sum_{q=0}^\infty \frac{1}
{\xi^{p+q}}
 \cdot\frac{2p+q}{p+q}\begin{pmatrix}p+q\cr q
\end{pmatrix}
\int D^p \cr&=&G_0(\xi)+\frac{1}{4\xi}\sum_{p=1}^\infty
\left(-\frac{s}{\xi}\right)^p \frac{2\xi-1}{(\xi-1)^{p+1}} \int
D^p\cr &=&G_0(\xi)+\frac{1-2\xi}{4(\xi-\xi^2)}\sum_{p=1}^\infty
\left(\frac{s}{\xi-\xi^2}\right)^p\int D^p
\end{eqnarray*}

We use the following formula for $G_0$, coming from Proposition 5.1:
$$G_0(\xi)=\frac{1}{2\xi}+\frac{1-2\xi}{4(\xi-\xi^2)}$$

Also, from $a^2+b^2+c^2+d^2=1$ we get $D=T-T^2$ with $T=a^2+b^2$, and we know from Lemma 5.1 that ${\rm law}(T)$ is the Lebesgue measure on
$[0,1]$, so:
\begin{eqnarray*}
\sum_{p=1}^\infty (-x)^p\int D^p
&=&\sum_{p=1}^\infty (-x)^p\sum_{k=0}^p
\begin{pmatrix}p\cr k\end{pmatrix}\int T^p(-T)^k\cr
&=&\sum_{p=1}^\infty (-x)^p\sum_{k=0}^p
\begin{pmatrix}p\cr k\end{pmatrix}\frac{(-1)^k}{p+k+1}\cr
&=&-1+\frac{4\,{\rm arcsinh}(\sqrt{x}/2)}{\sqrt{4x+x^2}}
\end{eqnarray*}

With $x=s/(\xi^2-\xi)$, we get the formula in the statement.
\end{proof}

We should mention that at $t=0,1$ the formula in Theorem 7.1 gives indeed  those in Proposition 5.1. At $t=0$ this follows from ${\rm arcsinh}(y)=\log(y+\sqrt{1+y^2})$, and at $t=1$ this follows from ${\rm arcsinh}(x)\sim x$ for $x\sim 0$.

The other remark is that Theorem 7.1 can be combined with the Stieltjes inverse formula, in order to get the density of the law of $w_t$. The corresponding function is piecewise analytic, and the precise formulae will not be given here.

\section{Poisson-Lebesgue interpolation}

We perform here a second computation, which is slightly more technical than the one in the previous section. Consider the following two elements:
\begin{eqnarray*}
v_0&=&\frac{1}{4}(u_{11}+u_{22}+u_{33}+u_{44})\cr
v_1&=&\frac{1}{2}(u_{11}+u_{22})
\end{eqnarray*}

We can consider the following element, interpolating between them:

\begin{definition}
To any real number $t$ we associate the element
$$v_t=\frac{1+t}{4}(u_{11}+u_{22})+\frac{1-t}{4}(u_{33}+u_{44})$$
where $u$ is the fundamental corepresentation of $A_s(4)$.
\end{definition}

The characteristic polynomial of $v_t$ can be computed by using
the matrix formulae in section 5. We get that this is the product
of the following polynomials:
$$Q_1(y)=y^2-(a^2+b^2)y+(1-t^2)a^2b^2$$
$$Q_2(y)=y^2-(c^2+d^2)y+(1-t^2)c^2d^2$$

Thus $v_t$ decomposes as a sum of two $2\times 2$ matrix
blocks, and we can compute its law, provided that we know
how to integrate polynomials in $a,b$.

We denote by $x=(a,b,c,d)$ the points on the real sphere $S^3$.

\begin{lemma}
We have the formula
$$\int a^{2k-2p}b^{2p}=\frac{4^{-k}}{(k+1)!}\cdot\frac{(2p)!(2k-2p)!}{p!(k-p)!}$$
where the integral is with respect to the uniform measure on
$S^3$.
\end{lemma}

\begin{proof}
We denote by $c_{rp}$ the numbers in the statement, where $r=k-p$:
$$c_{rp}=\int a^{2r}b^{2p}$$

Consider the following variables, depending on a real parameter
$t$:
\begin{eqnarray*}
A&=&\cos t\cdot a+\sin t\cdot b\cr B&=&-\sin t\cdot a+\cos t\cdot
b
\end{eqnarray*}

Since the map $(a,b)\to (A,B)$ is a rotation, we have:
$$c_{rp}=\int A^{2r}B^{2p}$$

We use the following formula, coming from $A'=B$ and $B'=-A$:
$$(A^{2r+1}B^{2p-1})'=(2r+1)A^{2r}B^{2p}-(2p-1)A^{2r+2}B^{2p-2}$$

Now since each $c_p$ doesn't depend on $t$, the integral of
$A^{2r+1}B^{2p-1}$ doesn't depend on $t$ either. Thus the
derivative of this integral vanishes, and we get:
$$(2r+1)\int A^{2r}B^{2p}=(2p-1)\int A^{2r+2}B^{2p-2}$$

This gives the following formula for the numbers $c_{rp}$:
\begin{eqnarray*}
c_{rp} &=&\frac{2p-1}{2r+1}\,c_{r+1,p-1}\cr
&=&\frac{(2p-1)(2p-3)\ldots (2p-2s+1)}{(2r+1)(2r+3)\ldots
(2r+2s-1)}\,c_{r+s,p-s}\cr &=&\frac{(2p-1)(2p-3)\ldots
1}{(2r+1)(2r+3)\ldots (2k-1)}\,c_{k0}\cr
\end{eqnarray*}

The number $c_{k0}$ is the $k$-th moment of $a$, so it is $4^{-k}$
times the $k$-th Catalan number, because $2a$ is known to be
semicircular. This gives the following formula:
\begin{eqnarray*}
c_{rp}&=&\frac{(2p)!}{2^p\,p!}\left(\frac{2^k\,k!}{(2k)!}\cdot\frac{(2r)!}{2^r\,r!}\right)\left(4^{-k}\cdot\frac{(2k)!}{k!(k+1)!}\right)\cr
&=&\frac{4^{-k}}{(k+1)!}\cdot\frac{2^k(2p)!(2r)!}{2^{p+r}p!r!}
\end{eqnarray*}

With $r=k-p$, this gives the formula in the statement.
\end{proof}

\begin{theorem}
The Cauchy transform of the law of $v_t$ is given by:
$$G'(\xi)=\frac{1}{2}\cdot\frac{2\xi-1}{\xi^2-\xi^3}
\sqrt{\frac{\xi-\xi^2}{\xi-\xi^2-(1-t^2)/4}}$$
\end{theorem}

\begin{proof}
By symmetry reasons the two blocks of $v_t$ have the
same spectral measure, so the Cauchy transform of $v_t$ is equal
to the Cauchy transform of each block, say of the first block. We
can apply Lemma 7.1, and we get:
\begin{eqnarray*}
G(\xi)&=&\frac{1}{\xi}+\frac{1}{2\xi}\sum_{p+q>0}\frac{(-1)^p}
{\xi^{2p+q}}\cdot\frac{2p+q}{p+q}\begin{pmatrix}p+q\cr q
\end{pmatrix}
\int (a^2+b^2)^q(sa^2b^2)^p\cr &=& \frac{1}{\xi}+\frac{1}{2\xi}
\sum_{p+q>0}\sum_{r=0}^q \frac{(-s)^p}{\xi^{2p+q}}
 \cdot\frac{2p+q}{p+q}\begin{pmatrix}p+q\cr q
\end{pmatrix}\begin{pmatrix}q\cr r
\end{pmatrix}
\int a^{2p+2q-2r}b^{2p+2r}\cr &=& \frac{1}{\xi}+\frac{1}{2\xi}
\sum_{p+q>0}\sum_{r=0}^q
\frac{(-s)^p}{\xi^{2p+q}}\cdot\frac{2p+q}{p+q}\,K
\end{eqnarray*}

Here $s=1-t^2$, and $K$ is the following product, obtained by
expanding all binomial coefficients, then by multiplying by the
quantity in Lemma 8.1:
\begin{eqnarray*}
K&=& \frac{(p+q)!}{p!q!}\cdot\frac{q!}{r!(q-r)!}\cdot
\frac{4^{-2p-q}}{(2p+q+1)!}\cdot\frac{(2p+2r)!(2p+2q-2r)!}
{(p+r)!(p+q-r)!}\cr
&=&\frac{4^{-2p-q}(p+q)!}{p!(2p+q+1)!}\cdot\frac{(2p+2r)!(2p+2q-2r)!}
{r!(q-r)!(p+r)! (p+q-r)!}
\end{eqnarray*}

We use now the following standard identity:
$$\sum_{r=0}^q\frac{(2p+2r)!(2p+2q-2r)!}
{r!(q-r)!(p+r)! (p+q-r)!}
=\frac{4^q}{q!}\cdot\frac{(2p)!(2p+q)!}{p!p!}$$

Thus when summing $K$ over $r$, we get the following quantity:
\begin{eqnarray*}
K_{r} &=&\frac{4^{-2p-q}(p+q)!}{p!(2p+q+1)!}
\cdot\frac{4^q}{q!}\cdot\frac{(2p)!(2p+q)!}{p!p!}\cr &=&
\frac{16^{-p}}{2p+q+1} \cdot \frac{(2p)!(p+q)!}{q!p!p!p!}
\end{eqnarray*}

We can get back now to the Cauchy transform:
\begin{eqnarray*}
G(\xi) &=& \frac{1}{\xi}+\frac{1}{2\xi} \sum_{p+q>0}
\frac{(-s)^p}{\xi^{2p+q}}\cdot\frac{2p+q}{p+q}\, K_r\cr &=&
\frac{1}{\xi}+\frac{1}{2\xi}\sum_{p+q>0}
\frac{(-s)^p}{\xi^{2p+q}}\cdot\frac{2p+q}{p+q}
\cdot\frac{16^{-p}}{2p+q+1} \cdot \frac{(2p)!(p+q)!}{q!p!p!p!} \cr
&=&\frac{1}{\xi}+\frac{1}{2}\sum_{p+q>0}
\left(-\frac{s}{16}\right)^p\frac{1}{\xi^{2p+q+1}}\cdot\frac{2p+q}{2p+q+1}
\cdot\frac{(2p)!(p+q-1)!}{q!p!p!p!}
\end{eqnarray*}

We can get rid of $p+q>0$ by using $G_0$, the value of $G$ at $s=0$:
\begin{eqnarray*}
G(\xi) &=&G_0(\xi)+\frac{1}{2}\sum_{p=1}^\infty \sum_{q=0}^\infty
\left(-\frac{s}{16}\right)^p\frac{1}{\xi^{2p+q+1}}\cdot\frac{2p+q}{2p+q+1}
\cdot\frac{(2p)!(p+q-1)!} {q!p!p!p!}\cr
&=&G_0(\xi)+\frac{1}{2}\sum_{p=1}^\infty\sum_{q=0}^\infty
\left(-\frac{s}{16}\right)^p\frac{(2p+q+1)^{-1}}{\xi^ {2p+q+1}}
\cdot \frac{2p+q}{p+q}
\begin{pmatrix}2p\cr
p\end{pmatrix}\begin{pmatrix}p+q\cr p\end{pmatrix}
\end{eqnarray*}

We take the derivative with respect to $\xi$:
\begin{eqnarray*}G'(\xi)&=&G_0'(\xi)-\frac{1}{2}\sum_{p=1}^\infty\sum_{q=0}^\infty
\left(-\frac{s}{16}\right)^p\frac{1}{\xi^{2p+q+2}}\cdot\frac{2p+q}{p+q}\begin{pmatrix}2p\cr
p\end{pmatrix}
\begin{pmatrix}p+q\cr p\end{pmatrix}\cr
&=&G_0'(\xi)-\frac{1}{2}\sum_{p=1}^\infty
\left(-\frac{s}{16}\right)^p\frac{1}{\xi^{p+2}}
\begin{pmatrix}2p\cr p\end{pmatrix}\sum_{q=0}^\infty
\frac{1}{\xi^{p+q}}\cdot\frac{2p+q}{p+q}
\begin{pmatrix}p+q\cr p\end{pmatrix}
\end{eqnarray*}

We use now the following standard identity:
$$\sum_{q=0}^\infty \frac{1}{\xi^{p+q}}\cdot\frac{2p+q}{p+q}\begin{pmatrix}p+q\cr
p\end{pmatrix}=\frac{2\xi-1}{(\xi-1)^{p+1}}$$

This gives the following formula for the Cauchy transform:
\begin{eqnarray*}
G'(\xi) &=&G_0'(\xi)-\frac{1}{2}\sum_{p=1}^\infty
\left(-\frac{s}{16}\right)^p\frac{1}{\xi^{p+2}}
\begin{pmatrix}2p\cr p\end{pmatrix}
\frac{2\xi-1}{(\xi-1)^{p+1}}\cr &=&
G_0'(\xi)+\frac{1}{2}\sum_{p=1}^\infty
\left(\frac{s}{16}\right)^p\frac{1}{\xi^{p+2}}
\begin{pmatrix}2p\cr p\end{pmatrix}
\frac{2\xi-1}{(1-\xi)^{p+1}}\cr &=&G_0'(\xi)+
\frac{2\xi-1}{2\xi^2-2\xi^3} \sum_{p=1}^\infty
\left(\frac{s/16}{\xi-\xi^2}\right)^p\begin{pmatrix}2p\cr
p\end{pmatrix}
\end{eqnarray*}

We use now the following standard identity:
$$\sum_{p=0}^\infty x^p\begin{pmatrix}2p\cr
p\end{pmatrix}=\frac{1}{\sqrt{1-4x}}$$

This gives the following formula for the Cauchy transform:
\begin{eqnarray*}
G'(\xi) &=&G_0'(\xi)+\frac{2\xi-1}{2\xi^2-2\xi^3}
\left(\left(1-4\cdot\frac{s/16}{\xi-\xi^2}\right)^{-1/2}-1
\right)\cr &=&G_0'(\xi)+\frac{2\xi-1}{2\xi^2-2\xi^3}
\left(\frac{\xi-\xi^2}{\xi-\xi^2-s/4}\right)^{1/2}-
\frac{2\xi-1}{2\xi^2-2\xi^3}
\end{eqnarray*}

We compute $G_0'$ by using Proposition 5.1, and we get:
$$G_0'(\xi)=-\frac{1}{2}\cdot\frac{2\xi-1}{\xi^3-\xi^2}$$

This gives the following formula for the Cauchy transform:
$$G'(\xi)=\frac{1}{2}\cdot\frac{2\xi-1}{\xi^2-\xi^3}
\sqrt{\frac{\xi-\xi^2}{\xi-\xi^2-s/4}}$$

Together with $s=1-t^2$, we get the formula in the statement.
\end{proof}

\section{Symmetric groups}

Let $u$ be the fundamental corepresentation of $C(S_4)$.
We consider the following element of $C(S_4)$, depending
on real parameters $t_i$ which sum up to $1$:
$$u_t=t_1u_{11}+t_2u_{22}+t_3u_{33}+t_4u_{44}$$

These can be regarded as being ``classical analogues'' of the variables considered in the previous sections. Their laws can be computed as follows.

\begin{theorem}
The law of $u_t$ is the following average of Dirac masses:
$${\rm law}(u_t)=\frac{1}{24}\left(9\delta_0+\delta_1+2\sum_i
\delta_{t_i}+\sum_{i\neq j} \delta_{t_i+t_j}\right)$$
\end{theorem}

\begin{proof}
We have $u_{ii}=\chi(X_i)$, where $X_i$ is the set of permutations
in $S_4$ fixing $i$. There are $6$ such permutations, namely:
\begin{enumerate}
\item  The identity $I$.

\item The two $3$-cycles fixing $i$; we denote by $C_i$ the set
they form.

\item The three transpositions fixing $i$; we denote by $T_i$ the
set they form.
\end{enumerate}

Observe that: the identity $I$ belongs to each $X_{i}$; the set
$C_i$ doesn't intersect the set $X_j$, for $j\neq i$; each of the
six transpositions in $S_4$ can be obtained by taking
intersections between the sets $T_i$ and their complements.

These remarks show that the algebra $\Delta$ generated by the
diagonal elements $u_{ii}$ is a $12$-dimensional vector space,
with the following basis:
\begin{enumerate}
\item The projection $\chi\{I\}$.

\item The $4$ projections $\chi\{C_i\}$.

\item The $6$ projections $\chi\{T_{ij}\}$, where $T_{ij}$ is the
transposition fixing $i\neq j$.

\item The projection $\chi\{ D\}$ onto what's left.
\end{enumerate}

With these notations, we have the following formula:
$$u_{ii}=\chi\{I\}+\chi\{C_i\}+\sum_{i\neq j}\chi\{T_{ij}\}$$

We get in this way a formula for average in the statement:
$$u_t=\chi\{I\}+\sum_it_i
\,\chi\{C_i\}+\sum_{i\neq j}(t_i+t_j)\,\chi\{T_{ij}\}$$

On the other hand, the restriction of the integration (or
averaging) over $S_4$ to the subalgebra $\Delta\subset \mathbb
C(S_4)$ is obtained by counting elements in various subsets of
$S_4$ corresponding to the above basis of $\Delta$. We get:
$$\int\chi\{I\}=\frac{1}{24},\int\chi\{C_i\}=\frac{1}{12},\int\chi(T_{ij})=\frac{1}{24},\int\chi\{D\}=\frac{9}{24}$$

We can compute now the moments of $u_t$:
\begin{eqnarray*}
\int u_t^k&=&\int \chi\{I\}+\sum_it_i^k \,\chi\{C_i\}+\sum_{i\neq
j}(t_i+t_j)^k\,\chi\{T_{ij}\}\cr
&=&\frac{1}{24}+\frac{1}{12}\sum_it_i^k +\frac{1}{24}\sum_{i\neq
j}(t_i+t_j)^k\cr &=&\frac{1}{24}\left( 1^k+2\sum_it_i^k
+\sum_{i\neq j}(t_i+t_j)^k+9\cdot 0^k\right)
\end{eqnarray*}

In this formula the $9\cdot 0^k=0$ term is there for the $9$
coefficient to produce the equality $24=1+2\cdot 4+6+9$. This
gives the formula in the statement.
\end{proof}

\begin{corollary}
The laws of the averages $m_s=(u_{11}+\ldots +u_{ss})/s$ are:
\begin{eqnarray*}
{\rm law}(m_1)&=&\frac{1}{24}\left(18\delta_0+6\delta_1\right)\cr
{\rm
law}(m_2)&=&\frac{1}{24}\left(14\delta_0+8\delta_{1/2}+2\delta_1\right)
 \cr {\rm
law}(m_3)&=&\frac{1}{24}\left(11\delta_0+9\delta_{1/3}+3\delta_{2/3}+
\delta_1\right) \cr {\rm
law}(m_4)&=&\frac{1}{24}\left(9\delta_0+8\delta_{1/4}+6\delta_{1/2}+
\delta_1\right)
\end{eqnarray*}
\end{corollary}

The challenging problem here is to work out the analogy with
$A_s(4)$. It is known from \cite{bc2} that the analogy between
classical and quantum appears in the limit $n\to\infty$, with the
Poisson semigroup of measures for $C(S_n)$ corresponding
to the free Poisson semigroup for $A_s(n)$. The above computations
should be regarded as a first step towards understanding what
happens with the analogy, when $n$ is fixed.

\end{document}